**Using the Blackwell-Cover Guessing Strategy to Separate Unobserved Bernoulli Trials - II**

By James D. Stein

Department of Mathematics and Statistics, California State University (Long Beach)

**Abstract -** In a prior paper, we investigated what could be said about the number in the pair that was not observed when we employed the Blackwell-Cover Guessing Strategy. In this paper, we generalize prior results to arbitrary distributions.

We use a variant of the Blackwell-Cover strategy ro take pairs of coins chosen from a mixture of two different distributions G and H such that one of the coins is flipped, and as a result of that flip the unobserved (unflipped) coin is placed in one of two sets A and B. We show that, regardless of which distributions are used for G and H, there are situations in which this can be done so that the mean of the heads probabilities of the unobserved coins in A is less than the mean of the heads probabilities of the unobserved coins in B.

What makes this result surprising, in addition to its generality, is a comparison with a corresponding result for marbles. If one knows that one will always be weighing one of two marbles with different specific weights, such as 10g and 20g, knowing the weight of the marble being weighed automatically determines the weight of the other. However, even if one knows that one will be flipping one of two coins with specific heads probabilities such as 1/3 and 2/3, the outcome of the flip of one of the coins will not immediately enable one to conclude anything about the heads probability of the unflipped coin – which is where the Blackwell-Cover Guessing Strategy comes into play.

**Section I - The Experimental Setup**

Short sections of this paper have been copied directly from [3] in order that this paper may be self-contained.

Let $0 \leq c \leq C \leq 1$ and $0 \leq d \leq D \leq 1$. In what follows, the individual conducting the experiment is simply presented with two identical-looking coins. One coin is labeled "Flip Me", and will be denoted X (the eXamined coin) and the other is unlabeled and will be denoted U (the Unexamined coin). The individual conducting the experiment flips the coin labeled "Flip Me". If it lands heads, he assigns the unexamined (unflipped) coin to set A, and if it lands tails, he assigns the unexamined coin to set B.

The individual or entity preparing the two coins to be given to the individual conducting the experiment has one bag of coins labeled G, which consists of coins with uniformly distributed heads probabilities on [c,C] and one bag of coins labeled H, which consists of coins with uniformly distributed heads probabilities on [d,D]. He also has a spinner which lands in the white region with probability $1/(1+\beta)$, where $\beta > 0$. It lands in the red region with probability $\beta/(1+\beta)$.

As will be seen, the computational results involve only m(G), the mean expected heads probability of a coin from G, m(H) (defined likewise) and $\beta$. We wish to construct an experiment in which the individual conducting that experiment cannot *a priori* conclude anything about the heads probability of an unflipped coin. If there were a way to distinguish coins coming from G from coins coming from H, the individual conducting the experiment might keep track of the results by distribution, and thus reach statistical conclusions about the means of the coins coming from G and H.

We initially present the results for uniform distributions G on [c,C] and H on [d,D] for two reasons. First, for anyone who is familiar with ([3]), the treatment parallels it and will make for easy reading. Second, because we can explain quite easily why the results obtained in this situation generalize in a straightforward fashion to arbitrary distributions G and H.

The individual or entity preparing the two coins follows one of two scenarios.

**Scenario I** – The spinner is spun twice. There are three possibilities.

(1) – Both spins land in the white region. In this case, two coins are chosen independently randomly from G and one is labeled "Flip Me".

(2) Both spins land in the red region. In this case, two coins are chosen independently randomly from H and one is labeled "Flip Me".

(3) One spin lands in the white region and one in the red region. In this case, one coin is chosen randomly from G and one coin is chosen randomly from H, the choice again being independent. A fair coin is flipped to determine which of these two coins is labeled "Flip Me".

**Scenario II** – The spinner is spun once. If it lands in the white region, a coin is chosen randomly (uniformly) from G and is labeled "Flip Me". The other coin is chosen randomly (uniformly) from H. If the spinner lands in the red region, the "Flip Me" coin is chosen randomly from H and the other coin from G.

In Scenario I, the pairs of coins occur with independent probabilities. The probability that the two coins are both from G is $1/(1+\beta)^2$, the probability that both coins are from H is $\beta^2/(1+\beta)^2$, and the probability that one coin is from G and the other from H is $2\beta/(1+\beta)^2$.

In Scenario II, every pair of coins consists of one each from both G and H. The ratio of coins from H that were flipped (examined) to coins from G that were flipped (examined) is β. The unseen coins therefore appear in the reciprocal ratio.

**An Intuitive Description of Scenarios I and II**

We think of the process of preparation of the two coins as a distribution constructed by thoroughly mixing 1 bag of coins with heads probabilities uniformly distributed in [c,C] (the distribution G) with β bags of coins with heads probabilities uniformly distributed in [d,D] (the distribution H). The coins are colored to denote which distribution, G or H, they come from (again, these colors are for the purpose of the reader or the preparer to distinguish between G and H and are not seen by the individual conducting the experiment). In both scenarios, first pick a coin from the bag with all the coins in it. In Scenario II, label this coin "Flip Me", and then pick the second coin, which is to be given to the conductor of the experiment, that is of the opposite color randomly from the bag. In Scenario I, pick another coin randomly from the bag. If the coins are of the same color, label a fair coin "Flip Me" and give it and one of the two coins to the conductor of the experiment. If the coins are of opposite color, flip a fair coin to determine which of the two coins is to be labeled "Flip Me" and give both coins to the conductor of the experiment.

**Section II – Computations for Scenario I**

We initially compute results for Scenario I, and then show how to modify these results to produce the corresponding formulas for Scenario II.

The expressions P(U→A) and m(U→A) denote the probability of assigning the unexamined coin to the set A and the mean of the heads probabilities of the coins assigned to A. The expression EH(A) denotes the expected number of heads assigned to A as a result of the flip. The expressions P(U→B), m(U→B) and EH(B) are defined similarly.

There are 8 separate cases corresponding to 3 binary choices – which distribution (G or H) the first coin to be flipped comes from, whether the flip is heads or tails, and which distribution (G or H) the second coin comes from. In the following table, **Dist 1** is the distribution (G or H) from which the first coin is taken, **Prob 1** is the probability that **Dist 1** is selected for the first coin, **Flip** is the outcome of the coin flip. **Dist 2** is the distribution (G or H) from which the second coin is taken, **Prob 2** is the probability that **Dist 2** is selected for the second coin, **Assign** is the outcome of the coin flip. q represents the heads probability for the first coin and p the heads probability for the second coin. **Comb Prob** is the combined probability that the unseen coin will be assigned to the set named in **Assign**, and **Exp Heads** is number of heads expected in the set named in **Assign** as the result of sending the unseen coin to that set. The integrations needed for the computations will be carried out following the presentation of the table.

**Table 1 (Values for a Single Coin)**

| Case | Dist 1 | Prob 1 | Flip | Dist 2 | Prob 2 | Assign | Comb Prob | Exp Heads |
|---|---|---|---|---|---|---|---|---|
| I | G | $1/(1+\beta)$ | heads | G | $1/(1+\beta)$ | A | $q/(1+\beta)^2$ | $qp/(1+\beta)^2$ |
| II | G | $1/(1+\beta)$ | tails | G | $1/(1+\beta)$ | B | $(1-q)/(1+\beta)^2$ | $(1-q)p/(1+\beta)^2$ |
| III | G | $1/(1+\beta)$ | heads | H | $\beta/(1+\beta)$ | A | $\beta q/(1+\beta)^2$ | $\beta qp/(1+\beta)^2$ |
| IV | G | $1/(1+\beta)$ | tails | H | $\beta/(1+\beta)$ | B | $\beta(1-q)/(1+\beta)^2$ | $\beta(1-q)p/(1+\beta)^2$ |
| V | H | $\beta/(1+\beta)$ | heads | G | $1/(1+\beta)$ | A | $\beta q/(1+\beta)^2$ | $\beta qp/(1+\beta)^2$ |
| VI | H | $\beta/(1+\beta)$ | tails | G | $1/(1+\beta)$ | B | $\beta(1-q)/(1+\beta)^2$ | $\beta(1-q)p/(1+\beta)^2$ |
| VII | H | $\beta/(1+\beta)$ | heads | H | $\beta/(1+\beta)$ | A | $\beta^2 q/(1+\beta)^2$ | $\beta^2 qp/(1+\beta)^2$ |
| VIII | H | $\beta/(1+\beta)$ | tails | H | $\beta/(1+\beta)$ | B | $\beta^2(1-q)/(1+\beta)^2$ | $\beta^2(1-q)p/(1+\beta)^2$ |

The integrations are straightforward – when a variable comes from G the limits of integration are c and C, and when the variable comes from H, the limits of integration are d and D. As both probabilities and expected heads are averages, the integrals are divided by the length of the appropriate integral as given by the standard formula for the average of a function on an interval, We show the integrals used to compute **Comb Prob** and **Exp Heads** for all coins assigned to B in case IV, and then simply list the results of the computation.

**Comb Prob** = $\frac{1}{(C-c)(D-d)(1+\theta)^2}\int_d^D \int_c^C \beta(1-q)dqdp$

**Exp Heads** = $\frac{1}{(C-c)(D-d)(1+\theta)^2}\int_d^D \int_c^C \beta p(1-q)dqdp$

In Table 1, **Assign** = A for cases I, III, V and VII. The computational results for these cases are obtained from evaluating integrals such as the ones above. m(G) and m(H) are the means of the uniform distributions on G and H respectively. It is their appearance in these formulas that enabled insights not apparent from the corresponding formulas in ([3]).

**Table 2 (Values for All Coins Assigned to A)**

| **Case** | **Comb Prob** | **Exp Heads** |
|---|---|---|
| I | $m(G)/(1+\beta)^2$ | $m(G)^2/(1+\beta)^2$ |
| III | $\beta m(G)/(1+\beta)^2$ | $\boldsymbol{\beta} m(G)m(H)/(1+\beta)^2$ |
| V | $\beta m(H)/(1+\beta)^2$ | $\boldsymbol{\beta} m(G)m(H)/(1+\beta)^2$ |
| VII | $\beta^2 m(H)/(1+\beta)^2$ | $\beta^2 m(H)^2/(1+\beta)^2$ |

The combined probability **(Comb Prob)** for each case can be thought of as the fraction of a coin assigned to A from that case. The expected heads (**Exp Heads**) can be thought of as the expected number of heads assigned to A from the fraction of the coin assigned to A (if 1/3 of a coin with heads probability ¼ is assigned to A, the expected number of heads from doing so is 1/3 x ¼ = 1/12).

P(U→A) can now be computed by adding up the fractional coins assigned to A in the four above cases. EH(A) can be computed by adding up the expected number of heads assigned to A in the four above cases. m(U→A) is then computed as EH(A)/P(U→A).

We obtain the following results.

$$P(U\to A) = \frac{m(G) + \beta m(H)}{1+\beta} \qquad EH(A) = \frac{(m(G) + \beta m(H))^2}{(1+\beta)^2} \qquad m(U\to A) = \frac{m(G) + \beta m(H)}{1+\beta} \qquad [1]$$

Notice that $m(U \to A)$ is the mean of the original distribution, where 1 copy of G is thoroughly mixed with $\beta$ copies of H. This result was also obtained in ([3]), and indicates that this particular method of employing the Blackwell-Cover technique for separation is ineffective in producing different means for A and B when G and H are mixed in such a way that the selection of pairs of coins is independent (as indicated by the **Prob 1** and **Prob 2** columns.

The above table for Cases II, IV, VI and VIII is only modified by replacing q by 1-q in the expressions in the **Comb Prob** and **Exp Heads** columns. Because the double integrals are simply products of single integrals, the result is that when the first integral is evaluated for the distribution G, the expression m(G) is replaced by 1-m(G); likewise for H. We therefore obtain

**Table 3 (Values for All Coins Assigned to B)**

| **Case** | **Comb Prob** | **Exp Heads** |
|---|---|---|
| II | $(1-m(G))/(1+\beta)^2$ | $m(G)(1-m(G))/(1+\beta)^2$ |
| IV | $\beta(1-m(G))/(1+\beta)^2$ | $\beta(1-m(G))m(H)/(1+\beta)^2$ |
| VI | $\beta(1-m(H))/(1+\beta)^2$ | $\beta m(G)(1-m(H))/(1+\beta)^2$ |
| VIII | $\beta^2 m(H)/(1+\beta)^2$ | $\beta^2 m(H)(1-m(H))/(1+\beta)^2$ |

Adding up the various cases, we obtain

$$P(U \to B) = \frac{(1-m(G))+\beta(1-m(H))}{1+\beta} \qquad EH(B) = \frac{(1-m(G))+\beta(1-m(H))(m(G)+\beta m(H))}{(1+\beta)^2} \qquad [2]$$

$$m(U \to B) = \frac{m(G)+\beta m(H)}{1+\beta} \qquad [3]$$

**Section III – Computations for Scenario II**

Although Scenario II consists only of cases III, IV, V and VI, we need to modify Table 1 to reflect the fact that once the first coin has been chosen from G (resp. H), the second coin must be chosen from H (resp. G).

**Table 4 (Values for a Single Coin)**

| Case | Dist 1 | Prob 1 | Flip | Dist 2 | Assign | Comb Prob | Exp Heads |
|---|---|---|---|---|---|---|---|
| III | G | 1/(1+β) | heads | H | A | q/(1+β) | qp/(1+β) |
| IV | G | 1/(1+β) | tails | H | B | (1-q)/(1+β) | (1-q)p/(1+β) |
| V | H | β/(1+β) | heads | G | A | β q/(1+β) | βqp/(1+β) |
| VI | H | β/(1+β) | tails | G | B | β(1-q)/(1+β) | β(1-q)p/(1+β) |

When we do the integrations for cases III and V (when coins are assigned to A), we obtain

**Table 5 (Values for All Coins Assigned to A)**

| Case | Comb Prob | Exp Heads |
|---|---|---|
| III | m(G)/(1+β) | m(G)m(H)/(1+β) |
| V | βm(H)/(1+β) | βm(G)m(H)/(1+β) |

These result in the following formulas for Scenario II (corresponding to Formula [1] for Scenario I).

$$P(U\to A) = \frac{m(G)+ \beta m(H)}{1+\beta} \qquad EH(A) = m(G)m(H) \qquad m(U\to A) = \frac{(1+\beta)m(G)m(H)}{m(G)+ \beta m(H)} \qquad [4]$$

We now do the integrations for cases IV and VI, when coins are assigned to B, we obtain

**Table 6 (Values for All Coins Assigned to B)**

| Case | Comb Prob | Exp Heads |
|---|---|---|
| IV | $(1\text{-}m(G))/(1+\beta)$ | $(1\text{-}m(G))m(H)/(1+\beta)$ |
| VI | $\beta(1\text{-}m(H))/(1+\beta)$ | $\beta m(G)(1\text{-}m(H))/(1+\beta)$ |

The results for Scenario II are expressed below.

$$P(U\rightarrow B) = \frac{(1-m(G))+\beta(1-m(H))}{1+\beta} \qquad EH(B) = \frac{(1-m(G))m(H)+\beta(1-m(H))m(G)}{1+\beta} \qquad [5]$$

$$m(U\rightarrow B) = \frac{m(H)+\beta m(G)-(1+\beta)m(G)m(H)}{(1-m(G))+\beta(1-m(H))} \qquad [6]$$

These results generalize those in ([4]), where the intervals G and H each consisted of a single point and $\beta = 1$.

**Theorem 1** - Let $m = (m(G) + \beta m(H))/(1+\beta)$, the mean of distribution D (the 1 to $\beta$ mixing of the distributions G and H). Then

1) $m(U\rightarrow A) \leq m(U\rightarrow B)$, with equality holding iff $m(G) = m(H)$.

2) $m(U\rightarrow A) < m$ iff $(m(G) - m(H))(m(G) - \beta^2 m(H)) > 0$.

3) $m(U\rightarrow B) > m$ iff $(m(G) - m(H))(\beta^2 (1- m(H)) - (1 - m(G))) > 0$.

**Proof:** As this is just algebra, let $x = m(G)$, $y= m(H)$.

1) Cross multiply; we then need to show that $(1+\beta)[(1-x)+\beta(1-y)] \leq (x+\beta y)(\beta x+y-(1+\beta)xy)$. Expand both sides and cancel $- (1+\beta)x^2y - \beta(1+\beta)xy^2$ from both sides. This reduces the problem to showing $(1+\beta)^2xy \leq \beta(x^2 + y^2) + (1 + \beta^2)xy$, which is equivalent to showing $0 \leq \beta(x-y)^2$.

2) Cross-multiplication reduces this to showing $(1+\beta)^2xy < (x + \beta y)^2$. Canceling $2\beta xy$ from both sides, regrouping and factoring yields the desired result.

Initially this may seem a somewhat surprising result; after all, the assignment algorithm was structured in an attempt to ensure that the unseen coin that most probably has the lower heads probability is assigned to A. However, if we assume $m(G) > m(H)$ and $\beta > 1$, we are more likely to examine a coin from H than a coin from G, and the algorithm "does the wrong thing" (see [3], p. 3) and assigns the coin more likely to have the higher heads probability to A.

3) Cross-multiply and cancel $\beta y - \beta x - 2\beta xy$ from both sides. Write what remains as inequality which is greater than 0 and factor. ■

**Corollary 1.1 –** (1) If $\beta = 1$, then $m(U\rightarrow A) < m < m(U\rightarrow B)$

(2) If $x = y$, then $m(U\rightarrow A) = m = m(U\rightarrow B)$

**Proof:** 1) Letting $x = m(G)$ and $y = m(H)$, both $m(U\rightarrow A) < m$ and $m < m(U\rightarrow B)$ reduce via cross-multiplication and simplification to $4xy < (x+y)^2$, which is equivalent to $(x-y)^2 > 0$. This is always true as long as $x \neq y$.

2) Substituting $x = m(G) = m(H)$ shows that $m(U\rightarrow A) = (1+\beta)x^2/((1+\beta)x) = x$, $m = (x+x)/2 = x$ and $m(U\rightarrow B) = (1+\beta)(x - x^2)/((1+\beta)(1 - x)) = x$. ■

Notice that Cor. 1.1 (2) implies that the Blackwell-Cover Guessing Strategy, as employed in this paper, cannot be used to separate independently chosen unseen members of an arbitrary distribution into subsets with different means. This does not foreclose on other implementations of the Blackwell-Cover Guessing Strategy achieving such a result, as there is a wide range of such strategies. These include flipping coins multiple times and competitive variations such as King-of-the-Hill to improve the accuracy of categorization of the flipped coin.

**Section IV – Modifying Scenario II to Create Different Means for A and B**

It can be seen from equations [1] and [4] that the difficulty in creating sets A and B with different means occurs when both coins from the pair are chosen from the same distribution. Intuitively, that is because Blackwell flipping does the wrong thing – assigns coins with lower heads probabilities to B and coins with higher heads probabilities to A. It is easy to see this in the case where both coins have the same probability X. If X < ½, then the coin that is flipped will be judged to have the higher heads probability with probability X, and so the unflipped coin (which also has heads probability X) will be assigned to A with probability X, and to B with probability 1-X; which is the opposite of what we are trying to accomplish.

One way to avoid this is to modify Scenario I, as indicated below.

**Scenario III** – The spinner is spun twice. There are three possibilities.

(1) – Both spins land in the white region. In this case, a fair coin is labeled "Flip Me" and the other coin is chosen randomly from G.

(2) Both spins land in the red region. In this case, a fair coin is labeled "Flip Me" and the other coin is chosen randomly from H.

(3) One spin lands in the white region and one in the red region. In this case, one coin is chosen randomly from G and one coin is chosen randomly from H. A fair coin is flipped to determine which of these two coins is labeled "Flip Me".

Scenario III causes the assigned coins from pairs chosen from the same distribution to be assigned equiprobably to A and B. This causes the first integration in the double integration to result in ½, rather than m(G) or m(H). This results in the following modifications to Tables 2 and 3.

**Modified Table 2 for Scenario III (Values for All Coins Assigned to A)**

| **Case** | **Comb Prob** | **Exp Heads** |
|---|---|---|
| I | $\frac{1}{2}/(1+\beta)^2$ | $\frac{1}{2}\,m(G)/(1+\beta)^2$ |
| III | $\beta m(G)/(1+\beta)^2$ | $\beta m(G)m(H)/(1+\beta)^2$ |
| V | $\beta m(H)/(1+\beta)^2$ | $\beta m(G)m(H)/(1+\beta)^2$ |
| VII | $\frac{1}{2}\,\beta^2/(1+\beta)^2$ | $\frac{1}{2}\,\beta^2 m(H)/(1+\beta)^2$ |

This results in the following changes to equation [1].

$$P(U\rightarrow A) = \frac{\frac{1}{2}(1+\beta^2)+\beta(m(G)+m(H))}{(1+\beta)^2} \qquad EH(A) = \frac{\frac{1}{2}\left(m(G)+\beta^2 m(H)\right)+2\beta m(G)\,m(H)}{(1+\beta)^2} \qquad [7]$$

$$m(U\rightarrow A) = \frac{(m(G)+\beta^2 m(H))+4\beta m(G)\,m(H)}{(1+\beta^2)+2\beta(m(G)+m(H))} \qquad [8]$$

**Modified Table 3 for Scenario III (Values for All Coins Assigned to B)**

| Case | Comb Prob | Exp Heads |
|---|---|---|
| II | $\frac{1}{2}/(1+\beta)^2$ | $\frac{1}{2}\,m(G)/(1+\beta)^2$ |
| IV | $\beta(1-m(G))/(1+\beta)^2$ | $\beta(1-m(G))m(H)/(1+\beta)^2$ |
| VI | $\beta(1-m(H))/(1+\beta)^2$ | $\beta m(G)(1-m(H))/(1+\beta)^2$ |
| VIII | $\frac{1}{2}\beta^2/(1+\beta)^2$ | $\frac{1}{2}\beta^2 m(H)/(1+\beta)^2$ |

This results in the following changes to equations [2].and [3]/

$$P(U\rightarrow B) = \frac{\frac{1}{2}(1+\beta^2)+\beta((1-m(G))+(1-m(H)))}{(1+\beta)^2} \qquad [9]$$

$$EH(B) = \frac{\frac{1}{2}\left(m(G)+\beta^2 m(H)\right)+\beta((1-m(G))m(H)+m(G)(1-m(H)))}{(1+\beta)^2} \qquad [10]$$

$$m(U\rightarrow B) = \frac{(m(G)+\beta^2 m(H))+2\beta\big((1-m(G))m(H)+m(G)(1-m(H))\big)}{(1+\beta^2)+2\beta((1-m(G))+(1-m(H)))} \qquad [11]$$

**The Case β = 1 for Equations [8] and [11]**

The relations between m(U→A), m and m(U→B) expressed in Theorem 1 are considerably more complicated in the case of Scenario III. However, it is easy to show that in the important case $\beta = 1$, $m(U\rightarrow A) \leq m \leq m(U\rightarrow B)$. Letting x = m(G) and y = m(H), we see that m(U→A) = $\frac{(m(G)+m(H))+4m(G)\,m(H)}{(1+\beta^2)+2(m(G)+m(H))} = \frac{x+y+4xy}{2+2(x+y)}$. Since m = (x+y)/2, by cross-multiplication we need merely show that $2(x+y) + 8xy \leq 2(x+y) + 2(x+y)^2$ to show that $m(U\rightarrow A) \leq m$. It is easily seen that this reduces to $0 \leq (x - y)^2$.

Since m(U→B) = $\frac{(m(G)+m(H))+2(m(H)(1-m(G))+m(G)(1-m(H)))}{2+2((1-m(G))+(1-m(H)))} = \frac{3(x+y)-4xy}{4-2(x+y)}$. by cross-multiplication we need merely show that 6(x+y) - 8xy ≥ 4(x+y) – 2(x+y)$^2$ to show that m ≤ m(U→B). This merely requires showing that (x+y)$^2$ + x + y ≥ 4xy, or (x – y)$^2$ + x + y ≥ 0.

**Section V – Generalizing Equations [1] – [11] to Arbitrary Distributions**

We observed in Section II that the computations for Case IV required the evaluations of the double integrals below.

**Comb Prob** = $\frac{1}{(C-c)(D-d)(1+\theta)^2}\int_d^D \int_c^C \beta(1-q)dqdp$

**Exp Heads** = $\frac{1}{(C-c)(D-d)(1+\theta)^2}\int_d^D \int_c^C \beta p(1-q)dqdp$

All integral evaluations reduce to evaluating expressions of the form $\frac{1}{b-a}\int_a^b 1\, dz$ or $\frac{1}{b-a}\int_a^b z\, dz$ , where z is either p or q. The first expression evaluates to 1, the second to the mean for the uniform distribution on the appropriate interval. However, if we assume that the distributions are arbitrary on [0,1] (rather than a uniform distribution on an interval), we can simply use the appropriate density function f(z) to see that we will either be evaluating $\int_0^1 f(z)\, dz$ or $\int_0^1 zf(z)\, dz$ instead of $\frac{1}{b-a}\int_a^b 1\, dz$ or $\frac{1}{b-a}\int_a^b z\, dz$ respectively, and these evaluate to 1 and the mean (expected value) of the appropriate distribution.

.